\magnification\magstep1
\font\titlefont=cmbx12

\font\ninerm=cmr9  \font\eightrm=cmr8  \font\sixrm=cmr6
\font\ninei=cmmi9  \font\eighti=cmmi8  \font\sixi=cmmi6
\font\ninesy=cmsy9 \font\eightsy=cmsy8 \font\sixsy=cmsy6
\font\ninebf=cmbx9 \font\eightbf=cmbx8 \font\sixbf=cmbx6
\font\nineit=cmti9 \font\eightit=cmti8 
\font\ninett=cmtt9 \font\eighttt=cmtt8 
\font\ninesl=cmsl9 \font\eightsl=cmsl8

\font\tengoth=eufm10  \font\ninegoth=eufm9
\font\eightgoth=eufm8 \font\sevengoth=eufm7 
\font\sixgoth=eufm6   \font\fivegoth=eufm5
\newfam\gothfam \def\goth{\fam\gothfam\tengoth} 
\textfont\gothfam=\tengoth
\scriptfont\gothfam=\sevengoth 
\scriptscriptfont\gothfam=\fivegoth

\catcode`@=11
\newskip\ttglue

\def\tenpoint{\def\rm{\fam0\tenrm}
  \textfont0=\tenrm \scriptfont0=\sevenrm
  \scriptscriptfont0\fiverm
  \textfont1=\teni \scriptfont1=\seveni
  \scriptscriptfont1\fivei 
  \textfont2=\tensy \scriptfont2=\sevensy
  \scriptscriptfont2\fivesy 
  \textfont3=\tenex \scriptfont3=\tenex
  \scriptscriptfont3\tenex 
  \textfont\itfam=\tenit\def\it{\fam\itfam\tenit}%
  \textfont\slfam=\tensl\def\sl{\fam\slfam\tensl}%
  \textfont\ttfam=\tentt\def\tt{\fam\ttfam\tentt}%
  \textfont\gothfam=\tengoth\scriptfont\gothfam=\sevengoth 
  \scriptscriptfont\gothfam=\fivegoth
  \def\goth{\fam\gothfam\tengoth}
  \textfont\bffam=\tenbf\scriptfont\bffam=\sevenbf
  \scriptscriptfont\bffam=\fivebf
  \def\bf{\fam\bffam\tenbf}%
  \tt\ttglue=.5em plus.25em minus.15em
  \normalbaselineskip=12pt \setbox\strutbox\hbox{\vrule
  height8.5pt depth3.5pt width0pt}%
  \let\big=\tenbig\normalbaselines\rm}

\def\ninepoint{\def\rm{\fam0\ninerm}
  \textfont0=\ninerm \scriptfont0=\sixrm
  \scriptscriptfont0\fiverm
  \textfont1=\ninei \scriptfont1=\sixi
  \scriptscriptfont1\fivei 
  \textfont2=\ninesy \scriptfont2=\sixsy
  \scriptscriptfont2\fivesy 
  \textfont3=\tenex \scriptfont3=\tenex
  \scriptscriptfont3\tenex 
  \textfont\itfam=\nineit\def\it{\fam\itfam\nineit}%
  \textfont\slfam=\ninesl\def\sl{\fam\slfam\ninesl}%
  \textfont\ttfam=\ninett\def\tt{\fam\ttfam\ninett}%
  \textfont\gothfam=\ninegoth\scriptfont\gothfam=\sixgoth 
  \scriptscriptfont\gothfam=\fivegoth
  \def\goth{\fam\gothfam\tengoth}
  \textfont\bffam=\ninebf\scriptfont\bffam=\sixbf
  \scriptscriptfont\bffam=\fivebf
  \def\bf{\fam\bffam\ninebf}%
  \tt\ttglue=.5em plus.25em minus.15em
  \normalbaselineskip=11pt \setbox\strutbox\hbox{\vrule
  height8pt depth3pt width0pt}%
  \let\big=\ninebig\normalbaselines\rm}

\def\eightpoint{\def\rm{\fam0\eightrm}
  \textfont0=\eightrm \scriptfont0=\sixrm
  \scriptscriptfont0\fiverm
  \textfont1=\eighti \scriptfont1=\sixi
  \scriptscriptfont1\fivei 
  \textfont2=\eightsy \scriptfont2=\sixsy
  \scriptscriptfont2\fivesy 
  \textfont3=\tenex \scriptfont3=\tenex
  \scriptscriptfont3\tenex 
  \textfont\itfam=\eightit\def\it{\fam\itfam\eightit}%
  \textfont\slfam=\eightsl\def\sl{\fam\slfam\eightsl}%
  \textfont\ttfam=\eighttt\def\tt{\fam\ttfam\eighttt}%
  \textfont\gothfam=\eightgoth\scriptfont\gothfam=\sixgoth 
  \scriptscriptfont\gothfam=\fivegoth
  \def\goth{\fam\gothfam\tengoth}
  \textfont\bffam=\eightbf\scriptfont\bffam=\sixbf
  \scriptscriptfont\bffam=\fivebf
  \def\bf{\fam\bffam\eightbf}%
  \tt\ttglue=.5em plus.25em minus.15em
  \normalbaselineskip=9pt \setbox\strutbox\hbox{\vrule
  height7pt depth2pt width0pt}%
  \let\big=\eightbig\normalbaselines\rm}
  
\newcount\newpen\newpen=50
\outer\def\beginsection#1\par{\vskip0pt
  plus.3\vsize\penalty\newpen\newpen=-50\vskip0pt
  plus-.3\vsize\bigskip\vskip\parskip
  \message{#1}\leftline{\bf#1}
  \nobreak\smallskip\noindent}

\def\QQ{{\bf Q}}
\def\RR{{\bf R}}
\def\ZZ{{\bf Z}}
\def\CC{{\bf C}}

\def\BK{{1}}
\def\BB{{2}}
\def\CH{{3}}
\def\CR{{4}}
\def\DAV{{5}}
\def\HA{{6}}
\def\HH{{7}}
\def\KK{{8}}
\def\LIG{{9}}
\def\JPS{{10}}
\def\SH{{11}}
\def\SIL{{12}}
\def\ST{{13}}
\def\STZ{{14}}
\def\STZZ{{15}}
\def\TZ{{16}}
\def\Z{{17}}

\noindent  \hfill June 23, 2011
\hrule\medskip\noindent
\bigskip\bigskip
\noindent \centerline{\titlefont  Integral points of a modular curve of level 11}
\bigskip
\noindent \centerline{\titlefont by  Ren\'e Schoof and Nikos Tzanakis}

\bigskip\bigskip\bigskip
\noindent{\bf Abstract.}  Using lower bounds for linear forms in elliptic logarithms we determine the  integral points of the modular curve  associated to the normalizer 
of a non-split Cartan group of level~$11$.
As an application we obtain a new solution of the class number one problem for complex quadratic fields.
\footnote{}{2000 Mathematics Subject Classification 11DXX, 11F03, 11G05, 11J86, 11R29} 

\beginsection 1. Introduction.

Let $E$ be the elliptic curve given by the Weierstrass equation
$$
Y^2+11Y\,=\,X^3+11X^2+33X.
$$
This is the curve 121B1 in Cremona's table~[\CR,~p.121]. By~[\CR,~p.256] the group of rational points of $E$ is the infinite cyclic group generated by the point~$(0,0)$. 
The main result of this paper is the following.

\proclaim Theorem 1.1. There are precisely seven rational points $(x,y)$ on $E$ for which
$${{x}\over{xy-11}}
$$
is integral. They are  $(0,0)$, $(0,-11)$, $(-2,-5)$, $(-2,-6)$, $(-6,-2)$, $(-11/4,-33/8)$ 
and the point at infinity. 

\medskip\noindent
Let $X_{ns}(11)$ denote the modular curve associated to the normalizer of a non-split Cartan subgroup of level~$11$; see~[\JPS, Appendix].
This curve has genus~$1$, is defined over $\QQ$ and parametrizes elliptic curves with a certain level~$11$ structure.
The interest of Theorem~1.1 lies in the fact that $X_{ns}(11)$ is isomorphic over $\QQ$ to the curve $E$ and that rational points $(x,y)$ on $E$
for which  $x/(xy-11)$ is integral, correspond  precisely to
{\it integral}  points on $X_{ns}(11)$, i.e., to rational points  for which the parametrized elliptic curve has its $j$-invariant in~$\ZZ$.
We prove this in section 4  and  hence obtain the following corollary of Theorem 1.1.

\proclaim Theorem 1.2. There are precisely seven integral points on the modular curve $X_{ns}(11)$. 

\medskip\noindent
As explained by J-P.~Serre in the Appendix of~[\JPS], every imaginary quadratic order $R$ of class number~1 in which the prime $11$ is inert, gives rise to an integral point on $X_{ns}(11)$. The elliptic curve that is parametrized by this point admits  complex multiplication (CM) by~$R$.
Since $11$ is inert in the quadratic orders of discriminant $-3$, $-4$, $-12$, $-16$, $-27$, $-67$ and $-163$, all of which have class number~1, the seven integral points of Theorem~1.1 are accounted for by these CM curves. See section~4 for the precise correspondence.

If the class number  of the imaginary quadratic order of discriminant $\Delta$ is~1, then the prime $11$ is inert in it whenever~$|\Delta|>44$. Therefore the fact that there are no other integral points on $X_{ns}(11)$ gives an independent proof of the Baker-Heegner-Stark theorem~[\BK, \HH, \ST]: the only imaginary quadratic orders with class number~1 are the ones with discriminant equal 
to one of $-3$,  $-4$,  $ -7$,  $ -8$,  $ -11$,  $ -12$,  $ -16$,  $-19$,  $ -27$,  $ -28$,  $ -43$,  $ -67$ and $-163$. 

Our proof exploits effective lower bounds for linear forms in elliptic logarithms~[\DAV].
In this respect it differs from earlier work by M.~Kenku~[\KK]  and B.~Baran~[\BB],
who exploit modular curves of level~$7$ and $9$ respectively.  In both  cases the curves involved have genus~0 and the  problem is reduced to a  cubic Thue equation, which is solved by Skolem's method.

The paper is organized as follows. In section 2  we prove two inequalities. These  are used in section~3, where we apply the method of linear forms in elliptic logarithms and prove Theorem~1.1. In section 4 we explain the relation with the modular curve of level~$11$ and prove Theorem 1.2.
All  calculations  can be checked easily and quickly  by means of the PARI software package.

\beginsection 2. Two inequalities.

In this section we prove two inequalities concerning the elliptic curve $E$ given by the Weierstrass equation $Y^2+11Y\,=\,X^3+11X^2+33X$.
The first inequality regards a property  of the group $E(\RR)$ of  real points of $E$, while the second  is concerned with heights of points in the group $E(\QQ)$ of rational points.

Let $t$ be the function on $E$ given by 
$$t=Y-{11\over X}.
$$ 
It has simple poles at the points $(0,0)$ and $(0,-11)$ and a pole of order~3 at infinity.
Its zero locus consists of five distinct points, which we call the {\it cusps} of $E$ because under the isomorphism  of section~4 they correspond to  the cusps of the modular curve $X_{ns}(11)$. The $x$-coordinates of the cusps are the zeroes of the polynomial
$X^5 + 11X^4 + 33X^3 - 121X - 121$. In particular, they are all real.
It follows that the cusps are contained in the group~$E(\RR)$.
The curve $E$ has only one connected component over $\RR$,
so that $E(\RR)$ is homeomorphic to a circle.

Writing $F(X,Y)=Y^2+11Y-X^3-11X^2-33X$, we define the function $g$ by
$$
g\,\,=\,\,\det\pmatrix{{{\partial t}\over{\partial X}}&
{{\partial t}\over{\partial Y}}\cr{{\partial F}\over{\partial X}}&{{\partial F}\over{\partial Y}}\cr}\,\,=\,\,3X^2+22X+33+{{11(2Y+11)}\over{X^2}}.
$$
It has poles of order 2 at $(0,0)$ and $(0,-11)$ and a pole of order~4 at infinity. It has eight zeroes on $E$, four of which are real.

\proclaim Lemma 2.1. Let $U$ be the subset of $E(\RR)$ given by
$$
U\,\,=\,\,\{P\in E(\RR): \hbox{$|t(P)|<{1\over{20}}$}\}.
$$
Then 
\item{(a)} the set $U$ is the disjoint union of five open intervals, each containing precisely one cusp;
\item{(b)} the function $g$
satisfies $|g(P)|\ge 1$ for every $P\in U$.

\smallskip\noindent{\bf Proof.}  In the proof all values of the functions $g$ and $t$ are given with an accuracy of two decimals only. 
The values of $t$ in the four real zeroes of $g$ on $E(\RR)$ are
equal to $-7.39$, $0.63$, $-0.16$ and $-23.06$ respectively.
  Since the absolute value of each of these numbers exceeds~${1\over{20}}$, the function $g$ has no zeroes in~$U$ and hence Lagrange's multiplier method ensures that $t$ assumes {\it no} extremal values in~$U$. 
It follows that $U$ is a union of open intervals~$I$, each of which contains {\it at most} one zero of~$t$. If $t$ were not to vanish on an interval~$I$, then its values on the boundary points of $I$ would either be both equal to $+{1\over{20}}$ or to~$-{1\over{20}}$. This is impossible as $t$ assumes no extremal values on~$I$. This shows that $t$ has {\it at least} one zero in $I$ and (a) follows.

To prove $(b)$, note that the values of $g$ in the five cusps are $9.75$, $-1.78$, $1.39$, $-3.79$ and $159.43$ respectively. Therefore $|g(P)|\ge 1$ for all points $P$ in a sufficiently small neigborhood of the cusps. We need to show that $U$ is such a neighborhood.
We saw in the proof of part (a) that   $g$ has no zeroes in $U$. Lagrange's multiplier method shows that  $g$ assumes its extremal values in  the zeroes of the function 
$$
\det\pmatrix{{{\partial g}\over{\partial x}}&
{{\partial g}\over{\partial y}}\cr{{\partial F}\over{\partial x}}&{{\partial F}\over{\partial y}}\cr}\,\,=\,\,(2Y+11)\left(6X+22-{{22(2Y+11)}\over{X^3}}\right)+{{22(3X^2+22X+33)}\over{X^2}}.
$$
This function has five zeroes in $E(\RR)$ and the  function $t$ assumes  the values $-3.60$, $0.34$, $-5.19$, $-0.44$ and $2.57$ in these zeroes. Since the absolute values of
these numbers  exceed ${1\over{20}}$, the zeroes  are not contained in~$U$ and  hence
$g$ assumes no extremal values on~$U$.

It follows that on each of the five intervals $I$ of part (a) the function $g$ is monotonous and assumes either only positive or only negative values.
This implies that on each $I$ we have $|g(P)|\ge {\rm min}(|g(z)|, |g(z')|)$ where $z,z'$ are the boundary points of~$I$. In our case the boundary points are given by the equation $t=Y-{{11}\over X}=\pm{1\over{20}}$. The values of $g$ in these points are given by $9.30$, $-2.05$, $1.63$, $-4.21$, $159.23$ (for the plus sign) and $10.18$, $-1.46$, $1.14$, $-3.39$ and $159.62$ (for the minus sign) respectively. The number with the smallest absolute value is $1.14$ which still exceeds~$1$. This proves the Lemma.

\bigskip
\noindent{\bf Remark.} The proof of Lemma~2.1 is related to the arguments in~[\TZ,~section 2] and [\STZZ, section 2.4]. In our case the situation  is relatively straightforward because all zeroes of the function $t$ are simple.

For any non-constant $f$ in the function field of $E$ and any point $P\in E(\QQ)$ we let
$$
H_f(P)\,\,=\,\,\prod_{p\le\infty}{\rm max}(1,|f(P)|_p),
$$
denote the {\it height} of $P$ with respect to~$f$. 
We let $h_f(P)=\log\,H_f(P)$ denote the {\it logarithmic height} of~$P$ with respect to the function~$f$.
The {\it canonical height} $\widehat{h}(P)$ of $P$ is defined as
${1\over{{\rm deg}\,f}}\lim_{n\rightarrow\infty}h_f(2^nP)/4^n$.  Here $f$ can be any even non-constant function on~$E$, for instance $f=X$.  See~[\SIL,~VIII].  The function  $t=Y-11/X$  is not even. We consider  the height function $h_t$  and compare it to the canonical height. 
For our purposes the following weak estimate is sufficient.

\proclaim Lemma 2.2. For every point $P\in E(\QQ)$ we have
$$
\widehat{h}(P)\,\,\le\,\,\hbox{${1\over 3}$}h_t(P)+4.52.
$$

\smallskip\noindent{\bf Proof.} Let $P=(x,y)\in E(\QQ)$. We first compare $h_t(P)$ to $h_X(P)$. For every finite prime $p$ we have
$$
\max(1,|x|_p)\,\,\le\,\,\max(1,|y-\hbox{${{11}\over x}$}|_p)^{2/3}.
$$
This is obvious when $|x|_p\le 1$. When $|x|_p>1$, the Weierstrass equation implies that ${\rm ord}_p(x)=-2k$ and ${\rm ord}_p(y)=-3k$ 
for some $k>0$. It follows that ${\rm ord}_p(y-{{11}\over x})=-3k$ and the inequality follows.

At the infinite prime we have
$$
\max(1,|x|)\,\,\le\,\,7\max(1,|y-\hbox{${{11}\over x}$}|)^{2/3}.
$$
This is obvious when $|x|\le 7$. If $|x|>7$, 
we observe  that $E$ has no real points with $x$-coordinate less than $-7$, so that we actually have $x>7$. Then we have 
$(y+{{11}\over 2})^2\ge (x^{3/2}+{{11}\over2})^2$ and hence
$|y-{{11}\over x}|\,\,\ge\,\,|y+{{11}\over 2}|-|{{11}\over 2}+{{11}\over x}|\,\,\ge\,\,x^{3/2}+{{11}\over2}-{{99}\over{14}}\,\,\ge\,\,({1\over 7}x)^{3/2}$.

Taking the product, it follows that $H_X(P)\le 7H_t(P)^{2/3}$ and hence 
$$h_X(P)\le \hbox{${2\over 3}$}h_t(P)+\log 7.
$$ 
To conclude the proof, we compare $h_X(P)$ to $\widehat{h}(P)$. Since the discriminant of $E$ is $11^3$ and its $j$-invariant is $-2^{15}$, Silverman's estimate~[\SH,~Thm.1.1] implies 
$\widehat{h}(P)\le{1\over 2}h_X(P)+3.54.$
Combining the two estimates gives
$$
\widehat{h}(P)\,\,\le\,\,\hbox{${1\over 2}$}\left(\hbox{${2\over 3}$}h_t(P)+\log 7\right)+3.54\,\,<\,\,\hbox{${1\over 3}$}h_t(P)+4.52,
$$
as required

\beginsection 3. The proof.

In this section we prove Theorem 1.1.  Our proof closely follows the strategy of~[\STZ].
Let $\omega$ denote the invariant differential ${{dX}\over{2Y+11}}$ of $E$. 
We define the {\it elliptic logarithm} of a point $P\in E(\RR)$ 
by
$$
\lambda(P)\,\,=\,\,\int_{\infty}^P\omega.
$$
Since $\lambda(P)$ depends on the path of integration in $E(\CC)$, it
is  only well defined up to the period lattice of~$E$. Since $P$ is in $E(\RR)$, there is a path of integration inside the real locus~$E(\RR)$.
Therefore $\lambda(P)$  is equal to  a real number modulo the period lattice and this real number is unique up to 
a multiple of the {\it real period}  
$$\Omega\,=\,\int_{c}^{\infty}{{dx}\over{\sqrt{q(x)}}}\,=\,4.8024\ldots
$$
Here $c=-6.8026\ldots$ denotes the unique real zero of  $q(x)=x^3+11x^2+33x+{{121}\over 4}$. It follows that $\lambda(P)$
is a well defined element of $\RR/\Omega\ZZ$. The map $P\mapsto\lambda(P)$ is a continuous group isomorphism $E(\RR)\longrightarrow\RR/\Omega\ZZ$. In order to avoid ambiguity, we assume that $\lambda(P)$ is a real number satisfying~$0\le \lambda(P)<\Omega$.

\proclaim Lemma 3.1.  For any cusp $Q$ of $E$ we have $\lambda(Q)={k\over{11}}\Omega$ for some integer~$k$.  

\smallskip\noindent{\bf Proof.} Any cusp $Q=(x,y)$ is contained in $E(\RR)$ so that $\lambda(Q)=r\Omega$ for some $r\in\RR$. Since we have $y=11/x$, the  $x$-coordinate of $Q$ is a zero of the polynomial $ p(X)\,\,=\,\, X^5 + 11X^4 + 33X^3 - 121X - 121$. One checks that
the $11$-division polynomial of $E$ is divisible by~$X^5 + 11X^4 + 33X^3 - 121X - 121$.  This implies that $\lambda(Q)={k\over{11}}\Omega$ for some $k\in\ZZ$, as required.

Alternatively, one can avoid the computation of the $11$-division polynomial and proceed as follows. The curve $E$ admits complex multiplication by the ring $\ZZ[{{1+\sqrt{-11}}\over 2}]$ and the kernel of the endomorphism $\sqrt{-11}$
is precisely the order~$11$ group $G$ generated by $\lambda^{-1}({1\over{11}}\Omega)$. Since the Galois group of $\overline{\QQ}$ over $\QQ$  preserves $G$, there is a unique monic degree 5 polynomial $q(X)\in \QQ[X]$ whose zeroes are precisely the $x$-coordinates of the points of $G$. By the Nagell-Lutz Theorem~[\SIL,~VII.3.4], each point $(x,y)\in G$ has the property that $11x$ is an algebraic integer.
Therefore we can compute $q(X)$ by calculating sufficiently accurate approximations to its roots. We find that $p(X)=q(X)$ and hence $Q\in G$. This proves the lemma.

Finally, the lemma also follows from the fact that Halberstadt's isomorphism~[\HA,~3.3] is known to map the cusps  of the modular curve $X_{ns}(11)$ to certain $11$-torsion points of~$E$. See section~4.
\bigskip

\noindent{\bf Proof of Theorem 1.1.} First we check that the only integers $k$ with $|k|\le 20$ for which there are points $P=(x,y)$ in $E(\QQ)$ with $x/(xy-11)$ equal to $k$ are $k=0$, $\pm 2$, $-6$ and~$-8$. These values of $k$ already account for the seven points listed in Theorem~1.1.  To prove the theorem, let $P=(x,y)$ be  a point in $E(\QQ)$ for which $|x/(xy-11)|$ is an integer exceeding~$20$. Since $t=Y-11/X$, we have~$|t(P)|<{1\over{20}}$.  By Lemma~2.1 there is a cusp $Q$ and an open interval $I\subset U$ containing both $P$ and~$Q$.

Let $\int_{Q}^P\omega$ denote the integral from $Q$ to $P$ of the invariant differential $\omega$
along a path inside the interval $I$. Then~$\int_{Q}^P\omega$ is real and we estimate its absolute value. Writing
$F(X,Y)=Y^2+11Y-X^3-11X^2-33X$, we have for every function $f$ on $E$ that 
$$
df\,\,=\,\,{\rm det}\!\pmatrix{{{\partial f}\over{\partial X}}&
{{\partial f}\over{\partial Y}}\cr{{\partial F}\over{\partial X}}&{{\partial F}\over{\partial Y}}\cr}\omega.
$$
In particular, taking $f$ equal to $t$, we find that $dt/\omega$ is equal to
the function $g=3X^2+22X+33+{{11(2Y+11)}\over{X^2}}$
of Lemma~2.1. Therefore 
$$
\int_{Q}^P\omega=\int_0^{t(P)}{{\omega}\over{dt}}dt=\int_0^{t(P)}{{dt}\over{g}}.
$$
By Lemma~2.1 we have  $|g(x,y)|\ge 1$ for all $(x,y)\in I$. Therefore we have $|\int_{Q}^P\omega|\le|\int_0^{t(P)}dt| =|t(P)|$. Since $\lambda(P)-\lambda(Q)\equiv \int_{Q}^P\omega$ modulo $\Omega\ZZ$,   there exists $n'\in\ZZ$ such that 
$$
|\lambda(P)-\lambda(Q)+n'\Omega|\,\le\, |t(P)|.
$$
By Lemma~3.1 we then have
$$
|n{{\Omega}\over{11}}-\lambda(P)|\le |t(P)|,\qquad\hbox{for some $n\in\ZZ$.}
$$
Since $1/t(P)$ is in $\ZZ$, we have that $h_t(P)=-\log|t(P)|$. Therefore Lemma~2.2 implies that
$$
|n{{\Omega}\over{11}}-\lambda(P)|< \exp({13.56-3\widehat{h}(P)}),\qquad\hbox{for some $n\in\ZZ$.}
$$
We write $P_0$ for the generator $(0,0)$ of the group $E(\QQ)$ so that $P=mP_0$ for some integer~$m$. Since $\widehat{h}(P)=m^2\widehat{h}(P_0)$ and $\widehat{h}(P_0)=0.04489\ldots$ this gives
$$
|n\Omega-m\lambda(11P_0)|\le 11\cdot \exp({13.56-0.13\cdot m^2}),\qquad\hbox{for some $n\in\ZZ$.}\eqno{(1)}
$$
On the other hand, since $P_0$ is not a torsion point, $n\Omega-m\lambda(11P_0)$ is a non-vanishing linear form in the elliptic logarithms $\Omega$ and $\lambda(11P_0)$.   We recall the explicit lower bound that Sinnou David obtained for such forms~[\DAV,~Th\'eor\`eme~2.1]. In David's notation
we have $K=\QQ$, $D=1$ and~$k=2$. The coefficients $\beta_0,\beta_1,\beta_2$ of his linear form are equal to $0,n,-m$ in our case. We have $u_1=\Omega$ with $\gamma_1$ equal to the point at infinity and $u_2=\lambda(11P_0)$ with $\gamma_2=11P_0$. It follows that David's constants $V_1$ and $V_2$ are given by $V_1=1.415\ldots\times 10^{27}$ and $V_2=7.98\ldots\times 10^{14}$.

David's estimates imply that when
$B={\rm max}\,(|m|,|n|)$ exceeds $V_1= 1.415\ldots\cdot 10^{27}$, then we have
$$
|n\Omega-m\lambda(11P_0)|\,\,>\,\,\exp\left(-7.658\times 10^{44}(\log B+1)(\log\log B +15\log(2)+1)^3\right).\eqno{(2)}
$$
We consider first the case  $|m|\ge 12$. Then the right hand side of inequality (1) is $<0.07$. Since $\Omega=4.8024\ldots$ 
 and $\lambda(11P_0)=3.5579\ldots$, this  easily implies that~$|m|\ge |n|$ and hence~$B=|m|$. We claim that 
$$
|m|\,\,<\,\,1.415\times 10^{27}.\eqno{(3)}
$$
Indeed, if the inequality is false, we may apply David's lower bound~(2).
Comparing the inequalities (1) and (2) one finds
$|m|\,\,<\,\,3.62\times10^{25}$, contradicting our assumption.

The bound on $|m|$ is very large. However, we can use it to obtain a better bound by observing that for $|m|\ge12$  the right hand side of $(1)$ is less than $0.4 {{\Omega}\over{|m|}}$. This leads to the inequality
$$
|{n\over m}-{{\lambda(11P_0)}\over{\Omega}}|\,\,<\,\,{{0.4}\over{m^2}}\,\,<\,\,{1\over{2m^2}},
$$
implying  that $n/m$ is a convergent $p_k/q_k$ of the continued fraction expansion of~$\lambda(11P_0)/\Omega$.  
By $(3)$ we must have that~$q_k<1.415\times 10^{27}$.
Using Zagier's algorithm~[\Z] we compute $\lambda(11P_0)$ and $\Omega$ with an accuracy of $60$ decimal digits and verify that for $k>55$ the convergents $p_k/q_k$  do not satisfy~$q_k<1.415\times 10^{27}$.  Note that replacing $\lambda(11P_0)/\Omega$ by its approximation $\sigma$ to $60$ decimal digits, does not affect the first $55$ convergents. This follows from 
the inequality 
$$|{n\over m}-\sigma|\le|{n\over m}-{{\lambda(11P_0)}\over{\Omega}}|+|{{\lambda(11P_0)}\over{\Omega}}-\sigma|<
{{0.4}\over{m^2}}+10^{-60}<{1\over{2m^2}}.
$$
On the other hand, one checks that for $k\le 55$ 
inequality $(1)$ does not hold when~$q_k\ge 12$. Indeed, one has
$$
|p_k\Omega-q_k\lambda(11P_0)|\,>\,11\cdot \exp({13.56-0.13\cdot q_k^2})
$$
for all $k\le 55$ for which $q_k\ge 12$. This contradicts our assumption that $|m|\ge 12$.

It remains to deal with the case $|m|< 12$. Inspection of the points $mP_0$ for $-12< m< 12$ shows that only for $m=-2$, $-1$, $0$, $1$, $2$, $3$ and $4$ the point $(x,y)=mP_0$ has the property that $x/(xy-11)$ is integral. In fact, these are the seven points $(-2,-6)$, $(0,-11)$, $\infty$, $(0,0)$, $(-2,-5)$, $(-11/4,-33/8)$ and $(-6,-2)$ respectively. So once again we  recover the seven points of Theorem~1.1.
This completes the proof of Theorem~1.1. 

\medskip
The continued fraction argument to reduce the upper bound for $|m|$ is particularly simple in our case because the rank of the Mordell-Weil group of $E$ is~1. In general, one employs a  lattice reduction algorithm that can handle lattices of higher rank. See~[\TZ] where the LLL algorithm is used.
\medskip\noindent{\bf Normalizations.} 
Our definition of the canonical height agrees with the one  given by Silverman~[\SIL,~VIII]. The canonical height used by the PARI and MAGMA programs is twice as large, while the canonical height used by Sinnou David~[\DAV] is three times ours. In a similar way, our definition of the real period $\Omega$ agrees with the one given by Silverman~[\SIL] and the one used by PARI. The one used by Zagier~[\Z] is twice as large.

\beginsection 4. The modular curve.

Let $X_{ns}(11)$ denote the modular curve associated to the normalizer of a non-split Cartan subgroup of level $11$.
It parametrizes elliptic 
curves with a certain level $11$ structure~[\JPS,~Appendix].
In 1977,  G.~Ligozat~[\LIG,~Proposition 4.3.8.1] showed that 
$X_{ns}(11)$  is isomorphic
to the genus 1 curve  given by the Weierstrass equation
$Y^2+Y=X^3-X^2-7X+10$. Replacing $X$ by $X+4$ and $Y$ by $Y+5$, we see that this curve is isomorphic to the curve $E$ given by
$$
Y^2+11Y\,\,=\,\,X^3+11X^2+33X.
$$
In this section we show that the $j$-invariant of an elliptic curve parametrized by a point $P=(x,y)\in E(\QQ)$ is in~$\ZZ$ if and only if
$x/(xy-11)\in\ZZ$. This shows that Theorem~1.2 follows from Theorem~1.1.

The curve $X_{ns}(11)$ admits a natural morphism $j:X_{ns}(11)\longrightarrow{\bf P}^1$, mapping a point $P$ of $X_{ns}(11)$ to the 
$j$-invariant  of the elliptic curve parametrized by $P$. The morphism $j$ has degree 55 and is defined over $\QQ$.  The formula for  the natural morphism from $E$ to the $j$-line depends on the choice of an isomorphism between  the modular curve $X_{ns}(11)$ and the elliptic curve $E$. Since translation by a rational point is a $\QQ$-rational automorphism of $E$, there are infinitely many such  choices.   We
follow Halberstadt~[\HA] and choose one of the two isomorphisms that map the five cusps of $X_{ns}(11)$ to the  zeroes of the degree 5 function~$t=Y-11/X$.  See [\HA,~section 3]. For formulas that are based on a different choice see~[\CH].  According to Halberstadt's  explicit formula, we have 
$$j(X,Y)={{h(X,Y)}\over{(XY-11)^{11}}},
$$ 
where $h(X,Y)$ is equal to
$$
(X^2 + 11X + 22)^3((11X^2 + 88X + 121)Y + 2X^4 + 55X^3 + 451X^2 + 1452X + 1452)^3g(X,Y)
$$
and $g(X,Y)$ is the polynomial
$$\eqalign{&
(6750X^8 + 337590X^7 + 5159935X^6 + 36807958X^5 + 145636931X^4 + 341425458X^3\cr
& + 474292533X^2 + 362189058X + 117523307)Y + 51975X^9 + 1746052X^8 + 24440064X^7\cr
& + 188870352X^6 + 892661770X^5 + 2692703508X^4 + 5217583888X^3 + 6299026712X^2\cr
& + 4320837279X + 1288408000.\cr}
$$
Our formula follows from Halberstadt's formula [\HA,~(2--1)]  by dividing his polynomial  $f_3^3f_4^3$ by~$f_5^2$. After replacing $X$ by $X+4$ and $Y$ by $Y+5$, the quotient is our polynomial~$g(X,Y)$. 

One checks that the seven points listed in Theorem~1.1 give rise to the $j$-invariants $-5280^3$, $66^3$,  $12^3$, $-3\cdot 160^3$, $-640320^3$,  $0$ and $2\cdot 30^3$ respectively.
These are precisely the $j$-invariants of the elliptic curves with complex multiplication by  the quadratic orders of discriminant $-67$, $-16$, $-4$, $-27$, $-163$, $-3$ and $-12$ respectively.

\medskip
\proclaim Theorem 4.1.  Let $(x,y)$ be a rational point on the elliptic curve $E$ given by the Weierstrass equation $Y^2+11Y=X^3+11X^2+33X$.  Then $j(x,y)$ is in $\ZZ$ if and only if $x/(xy-11)$ is in $\ZZ$.

\smallskip\noindent{\bf Proof.} We study integrality of $j(x,y)$ and $x/(xy-11)$ one prime $l$ at a time. It follows from the Weierstrass equation that we can write $x=r/t^2$ and $y=s/t^3$ for certain $r,s,t\in\ZZ$
satisfying ${\rm gcd}(rs,t)=1$. The denominators of both $j(x,y)=h(x,y)/(xy-11)^{11}$ and $x/(xy-11)$ divide a power of $rs-11t^5$. Therefore, if $l$ is a prime not dividing $rs-11t^5$, both $j(x,y)$ and $x/(xy-11)$ are integral at~$l$.
Let therefore $l$ be a prime that divides $rs-11t^5$. If $l$ divides $t$, then it divides $rs$, which it cannot. So $l$ does not divide~$t$. This implies that both $x$ and $y$ are $l$-integral and $l$ divides~$xy-11$.

Suppose $l\not=11$. Then $l$ does not divide $x$, so that 
$x/(xy-11)$ is {\it not} integral at~$l$. Substituting $Y=11/X$ in the Weierstrass equation we find that $l$ divides 
$p(x)$ where $p(X)=X^5 + 11X^4 + 33X^3 -121X - 121$.
Suppose now that $j(x,y)=h(x,y)/(xy-11)^{11}$ is integral at~$l$. Then $l$ divides $h(x,y)$. Substituting $Y=11/X$ in the polynomial $h(X,Y)$ 
and multiplying by $X^4$, we find that $l$ divides 
$r(x)$ where $r(X)$ is a certain degree 31 polynomial  in $X$ with integral coefficients.  Therefore $l$ divides the resultant of  $p(X)$ and $r(X)$, which one checks to be equal 
to $11^{63}$. This shows that $l=11$. This contradicts our assumption and we conclude that $j(x,y)$ is not integral at~$l$.

Finally, suppose $l=11$. Since $l$ divides $xy-11$, it also divides $xy$ and it follows from the Weierstrass equation that $11$ 
actually divides both $x$ and $y$. 
It follows that $11$ divides $xy-11$ exactly once so that $x/(xy-11)$ is {\it integral} at~$11$. To see that $j(x,y)=h(x,y)/(xy-11)^{11}$ is also integral at~$11$, we observe that the  exact power of $11$ dividing $(xy-11)^{11}$
is $11^{11}$. On the other hand, one checks that when both $x$ and $y$ are divisible by~$11$, the numerator  $h(x,y)$ is divisible by~$11^{14}$. 
Therefore $h(x,y)/(xy-11)^{11}$ is divisible by $11^3$ and hence $j(x,y)$ is certainly integral at~$11$.

This proves the Theorem.

\bigskip\bigskip
\def\bibliography#1\par{\vskip0pt
  plus.3\vsize\penalty-250\vskip0pt
  plus-.3\vsize\bigskip\vskip\parskip
  \message{Bibliography}\leftline{\bf
  Bibliography}\nobreak\smallskip\noindent
  \ninepoint\frenchspacing#1}

\bibliography
{\item{[\BK]} Baker, A.:  A remark on the class number of quadratic fields   {\sl Bull. London Math. Soc.}  {\bf 1}   (1969),   98--102.
\item{[\BB]} Baran, B.: A modular curve of level 9 and the class number one problem, {\sl Journal of
Number Theory} {\bf 129} (2009), 715--728.
\item{[\CH]} Chen, I. and Cummins, C.: Elliptic curves with non-split mod 11 representations,
{\sl Math Comp.} {\bf 73} (2004), 869--880.
\item{[\CR]} Cremona, J.: {\sl Algorithms for modular elliptic curves},
2nd ed, Cambridge University Press, Cambridge 1997.
\item{[\DAV]} David, S.: Minorations de formes lin\'eaires de logarithmes elliptiques, 
{\sl  M\'emoires Soc. Math. France} (N.S) {\bf 62} (1995).
\item{[\HA]} Halberstadt, E.: Sur la courbe modulaire $X_{\scriptstyle \rm ndep}(11)$,
{\sl Experimental Math.}, {\bf 7} (1998), 163--174.
\item{[\HH]} Heegner, K.: Diophantische Analysis und Modulfunktionen, {\sl Math. Zeit.} {\bf 59} (1952), 227--253.
\item{[\KK]} Kenku, M.A.: A note on the integral points of a modular curve of level 7, {\sl Mathematika} {\bf 32} (1985), 45--48.
\item{[\LIG]} Ligozat, G.: {\sl Courbes modulaires de niveau $11$}, in J-P. Serre
and D.B. Zagier Eds, {\sl Modular Functions of one variable V}, LNM {\bf 601}, 149--237, Springer-Verlag, 1977
\item{[\JPS]} Serre, J-P.: {\sl Lectures on the Mordell-Weil Theorem}, Aspects of Mathematics {\bf 15}, Vieweg, Braunschweig 1997.
\item{[\SH]} Silverman, J.: The difference between the Weil height and the canonical height
on elliptic curves, {\sl Math.~Comp.} {\bf 55} (1990), 723--743.
\item{[\SIL]} Silverman, J.: {\sl The arithmetic of elliptic curves}, 2nd ed, Graduate Texts in Mathematics {\bf 106}, Springer-Verlag, New York 2009.
\item{[\ST]} Stark, H.M.: On complex quadratic fields with class number equal to one, {\sl Trans. Amer. Math. Soc.} {\bf 122} (1966), 112--119.
\item{[\STZ]} Stroeker, R.J. and Tzanakis, N.:  Solving elliptic diophantine equations
by estimating linear forms in elliptic logarithms, {\sl Acta Arith.} {\bf 67} (1994), 177--196.
\item{[\STZZ]} Stroeker, R.J. and Tzanakis, N.: Computing all integer solutions of a genus 1 equation,
{\sl Math.~Comp.} {\bf 72}  (2003) 1917--1933.
\item{[\TZ]} Tzanakis, N.:  Solving elliptic diophantine equations by estimating linear
forms in elliptic logarithms. The case of quartic equations, {\sl Acta Arith.} {\bf 75} (1996), 165--190.
\item{[\Z]} Zagier, D.:  Large integral points on elliptic curves, {\sl Math. Comp.} {\bf 48} (1987), 425--436.}

\bigskip\bigskip
{\eightpoint\hbox{\hskip 3cm\vbox{\hsize 7cm
\noindent Ren\'e Schoof\par
\noindent Universit\`a di Rome ``Tor Vergata"\par
\noindent Dipartimento di Matematica\par
\noindent I-00133 Roma, Italy\par
\noindent schoof@mat.uniroma2.it}\hfill
\vbox{\hsize 7cm
\noindent  Nikos Tzanakis\par
\noindent  Department of Mathematics\par
\noindent  University of Crete\par
\noindent  Iraklion, Greece\par
\noindent  tzanakis@math.uoc.gr}}

\bye